\documentclass{amsart}

\usepackage{amssymb}

\theoremstyle{plain}

\numberwithin{equation}{section}

\usepackage[all, 2cell, dvips]{xy}
\usepackage{latexsym,amsmath, amsfonts, graphics, epsf, epic}
\usepackage{amsthm}
\usepackage{amssymb}
\usepackage{amsopn}
\usepackage{amscd}
\usepackage{color}
\usepackage{amsmath, amscd, amssymb, amsthm, xspace}

\theoremstyle{plain}

\theoremstyle{definition}

\theoremstyle{plain}
\newtheorem{thm}{Theorem}[section]
\newtheorem{cor}[thm]{Corollary}
\newtheorem{lem}[thm]{Lemma}
\newtheorem{prop}[thm]{Proposition}

\theoremstyle{definition}

\newtheorem{remark}[thm]{Remark}

\newtheorem{defn}[thm]{Definition}

\newtheorem{question}[thm]{Question}

\usepackage{amssymb}

\def\al{\alpha}

\def\l.l.o.{\it l.l.o}

\def\chiup{\raise 2pt\hbox{$\chi$}}

\DeclareMathOperator{\spn}{span}

\newcommand{\N}{\mathbb{N}}

\newcommand{\gr}{\operatorname{gr}}

\newcommand{\s}{p}                      
\newcommand{\sps}{P}                    

\newcommand{\st}{\text{ }|\text{ }}

\newcommand{\dspfrac}[2]{\displaystyle{\frac{#1}{#2}}}

\title[Filtered Algebraic Algebras] {Filtered Algebraic Algebras}
\author[Alon Regev]{Alon Regev}
\address{ Department of Mathematical Sciences, Northern Illinois
University, Watson Hall 320 DeKalb, Illinois 60115, USA}
\email{regev@math.niu.edu }

\date{\today}

\begin{document}

\maketitle

\begin{abstract}
Small and Zelmanov posed the question whether every element of a graded algebra over an uncountable field must be nilpotent, provided that the homogeneous elements are nilpotent.
This question has recently been answered in the negative by A. Smoktunowicz. In this paper we prove that the answer is affirmative for associated graded algebras of filtered algebraic algebras.
Our result is based on Amitsur's theorems on algebas over infinite fields.
\end{abstract}

\section{Introduction}

Throughout this paper, $k$ is a field and $A$ is a $k$-algebra.
\begin{defn}

\begin{enumerate}
\item If $A$ is a $k$-algebra, we say that $A$ is \emph{nil} if every element of $A$ is nilpotent.
\item If $A$ is a graded $k$-algebra, we say that $A$ is graded-nil if every homogeneous element of $A$ is nilpotent.
\end{enumerate}
\end{defn}

Clearly, every nil, graded algebra is graded-nil. The question of
whether the converse to this statement is true has been of some
interest. Bartholdi \cite{Bartholdi} constructed a graded-nil but
not nil algebra over a countable field. Small and Zelmanov \cite{SZ}
posed the question whether, for an uncountable field, every
graded-nil algebra is nil. This was answered in the negative by A.
Smoktunowicz \cite{Agata-g}, who constructed an algebra over any
field which is graded-nil but not nil.

Let $A$ be a filtered $k$-algebra. That is, there exist subspaces $F_0\subseteq F_1 \subseteq \ldots $ of $A$, such that $\bigcup_{i\ge 0} F_i=A$ and $F_i F_j \subseteq F_{i+j}$ for all $i, j\ge 0$.
The associated graded algebra of $A$ is defined as follows.
\begin{defn}\label{Def-gr}
If $A$ is a $k$-algebra filtered by $\{F_n\}_{n\ge0}$, then the \emph{associated graded algebra} of $A$ (with respect to
this filtration) is 
\begin{equation}
\gr(A)=F_0 \oplus \frac{F_1}{F_0} \oplus \frac{F_2}{F_1} \oplus
\ldots,
\end{equation}
with multiplication defined as follows. For $p,q\ge 0$, if
$a_p+F_{p-1}$ and\\ $a_q + F_{q-1}$, with $a_p \in F_p$ and $a_q\in
F_q$, are arbitrary homogeneous elements of\\ $\gr(A)$ (taking
$F_{-1}=\{0\}$) then
\begin{equation} \label{assoc-mult}
(a_p+F_{p-1})(a_q + F_{q-1}) = a_p a_q + F_{p+q-1}.
\end{equation}
This multiplication is extended linearly to arbitrary elements of $\gr(A)$.
\end{defn}

We define
\begin{equation}
(\gr A)_{\ge 1} =  \frac{F_1}{F_0}\oplus \frac{F_2}{F_1}\oplus \ldots
\end{equation}
Our main result states that

\begin{thm}\label{thetheorem}
If $k$ is uncountable and $A$ is algebraic then $(\gr A)_{\ge 1}$ is nil.
\end{thm}

Thus, being graded-nil does imply being nil, for
 associated graded algebras of algebraic algebras over uncountable fields. This result is based on the Amitsur's Theorem \ref{LBI-LBD.intro} below.

\begin{defn}
\begin{enumerate}
\item
If $A$ is a nil $k$-algebra, we say that $A$ is \emph{locally of
bounded index (LBI)} over $k$ if the elements of every finite-dimensional $k$-subspace of $A$ have bounded index of nilpotence.
\item
If $A$ is an algebraic $k$-algebra, we say that $A$ is { \it
locally of bounded degree (LBD)} over $k$ if the elements of every
finite-dimensional $k$-subspace of $A$ have bounded degree of
algebraicity.
\end{enumerate}
\end{defn}
In~\cite{AIF} Amitsur proved
\begin{thm}[{\cite[Theorem 5]{AIF}}]\label{LBI-LBD.intro}
Let $k$ be an uncountable field and let $A$ be a $k$-algebra. If
$A$ is nil then it is LBI. If $A$ is algebraic then it is $LBD$.
\end{thm}
In fact, Amitsur proved the following stronger statements:
\begin{thm}[{\cite[Corollary 7]{AIF}}]\label{LBI-LBD-stronger.intro}
Let $k$ be an uncountable field and let $A$ be a $k$-algebra. Then
any nil subspace of $A$ has bounded index, and any
algebraic subspace of $A$ has bounded degree.
\end{thm}
In \cite{}, we use the properties of order-symmetric polynomials to give an alternative proof of Theorem \ref{LBI-LBD-stronger.intro}. In this paper, we use these properties, along with Theorem \ref{LBI-LBD-stronger.intro} itself, to prove Theorem \ref{thetheorem}.
\\

This paper is part of the author's Ph.D. dissertation at the University of California, San Diego under the supervision of Lance Small, to whom he wishes to  express his warm thanks.

\section {The polynomials $p_{i_1,\ldots,i_m}$}

Throughout this paper $k$ is a field and $k\langle
x_1,\ldots,x_m\rangle$ is the associative algebra of the
noncommutative polynomials in $x_1,\ldots,x_m$. Note that if $A=k\langle
x_1,\ldots,x_m\rangle$ then
\begin{equation}
A=\bigoplus_{n=0}^\infty A_n,
\end{equation}
where $A_n$ are the homogeneous polynomials in $x_1,\ldots,x_m$ of
total degree $n$.

Given a sequence $i_1,\ldots,i_m$ of nonnegative integers, there
are $\frac{(i_1+ \ldots+i_m)!}{i_1!\cdots i_m!}$ different
noncommutative monomials in $x_1,\ldots,x_m$ of degree $i_j$ in
$x_j$, $1\le j\le m$. For example let $m=2=i_1=i_2$, then the
$\frac{(2+2)!}{2!2!}=6$ monomials of degree 2 in $x_1$ and in
$x_2$ are $x_1^2x_2^2$, \ $x_1x_2x_1x_2$, \ $x_1x_2^2x_1$, \ $x_2x_1^2x_2$, \
$x_2x_1x_2x_1$ and $x_2^2x_1^2$.

\medskip
The polynomials $p_{i_1,\ldots,i_m}(x_1,\ldots,x_m)$ are defined,
for any $i_1,\ldots,i_k \in \N$, as follows.
\begin{defn}\label{thedef}
We define $\s_{i_1,\ldots,i_m}(x_1,\ldots,x_m)$ to be the sum of
all the $\frac{(i_1+ \ldots+i_m)!}{i_1!\cdots i_m!}$ different
monomials consisting of exactly $i_j$ occurrences of $x_j$ for each
$j$. We take $p_{0,\ldots,0}(x_1,\ldots,x_m)=1$.
\end{defn}

\begin{defn} \label{DefSn}
For any $n\geq 0$, any field $k$ and any $x_1,\ldots,x_m$ denote
\begin{equation}
\sps_n(x_1,\ldots,x_m)=\spn_k\{\s_{i_1,\ldots,i_m}(x_1,\ldots,x_m)\st
i_1+ \ldots+i_m=n\}.
\end{equation}
\end{defn}

\medskip
The polynomials $p_{i_1,\ldots,i_m}$ can be used as follows. It is
well known that for \\$f=\alpha_1x_1+ \ldots+\alpha_m x_m$,
\begin{equation}\label{nc-binom}
f^n=\sum_{1\le j_1,\ldots,j_n\le m}(\al_{j_1}x_{j_1})\cdots
(\al_{j_n}x_{j_n}) =\sum_{1\le j_1,\ldots,j_n\le
m}\al_{j_1}\cdots\al_{j_n}\cdot x_{j_1}\cdots x_{j_n}
\end{equation}

\begin{lem}\label{invar3}
Let $\al_1,\ldots,\al_m\in k$ and let $f=\al_1x_1+ \ldots
+\al_mx_m$. Then
\begin{equation}
f^{n}=\sum_{i_1+ \ldots+i_m=n}  \al_1^{i_1}\cdots \al_m^{i_m}\cdot
p_{i_1,\ldots,i_m}(x_1,\ldots,x_m).
\end{equation}

\begin{proof}
This follows from \eqref{nc-binom} and Definition \ref{thedef},
since for each monomial $x_{j_1}\cdots x_{j_n}$ that appears in $p_{i_1,\ldots, i_m}$, the coefficient $\al_{j_1}\cdots\al_{j_n}$ in \eqref{nc-binom} is equal to $\al_1^{i_1}\cdots \al_m^{i_m}$.
\end{proof}
\end{lem}

Let $W=\spn_k\{(\al_1x_1+ \ldots+\al_mx_m)^n\mid
\al_1,\ldots,\al_m\in k \}$. By Lemma~\ref{invar3}, $W\subseteq
P_n(x_1,\ldots,x_m )$. We will show:
\begin{prop}\label{v.d.m.1}
If $|k|\ge n+1$ then $W=P_n(x_1,\ldots,x_m)$.
\end{prop}
This is a  consequence of the fundamental
Proposition~\ref{main.vandermonde.1} below. We include its proof
here for completion, although it has been used in the past (see \cite{Kap}, \cite{Schur}). We begin with the following generalization of a
basic ``Vandermonde argument".
\begin{lem}\label{v.d.m.2}
Let $V$ be a vector space over a field $k$, let $v_0,\ldots,v_d\in
V$ and let $\xi_0,\ldots,\xi_d \in k$ be $d+1$ distinct field
elements. Let $W$ be a subspace of $V$. If
\begin{equation} \label{VdmEq}
\sum^d_{i=0}\xi_j^i v_i\in W
\end{equation}
for all $j \in \{0,\ldots,d\}$, then $v_i\in W$ for all $i\in
\{0,\ldots,d\}$.
\begin{proof}
We write the condition \eqref{VdmEq} in matrix form as follows:
\begin{displaymath}
\left(
\begin{array}{cccc}
1  & \xi_0 & \cdots & \xi
_0^d \\
1  & \xi_1 & \cdots & \xi_1^d\\
1  & \vdots & \ddots& \vdots \\
1  & \xi_d & \cdots & \xi_d^d\\
\end{array}
\right) \left(
\begin{array}{c}
v_0 \\
v_1 \\
\vdots  \\
v_d  \\
\end{array}
\right)= \left(
\begin{array}{c}
w_0 \\
w_1 \\
\vdots  \\
w_d  \\
\end{array}
\right),
\end{displaymath}
where $w_j\in W$. The multiplying matrix is a Vandermonde matrix,
known to be invertible when the $\xi_i$ are distinct. Multiplying
by the inverse matrix, we obtain $v_i\in W$ for all $i\in
\{0,\ldots,d\}$.
\end{proof}
\end{lem}
\begin{prop}\label{main.vandermonde.1}
Let $k$ be a field with $|k|\ge n+1$, let $W\subseteq V$ be
$k$-vector spaces, and let
 $\{ w_{\mu_1,\ldots, \mu_m} \mid \mu_1+ \ldots+\mu_m =n \}$ be a set of vectors in
$V$.
 Assume that
\begin{equation}\label{multi.v.d.m.1}
\sum_{\mu_1+ \ldots +\mu_m =n}\al_1^{\mu_1}\cdots
\al_m^{\mu_m}w_{\mu_1,\ldots,\mu_m} \in W
\end{equation}
for all $\al_i\in k$. Then all $w_{\mu_1,\ldots,\mu_m} \in W$.
\end{prop}

\begin{proof}
We prove the claim by induction on $m$. When $m=1$, equation
\eqref{multi.v.d.m.1} becomes\\ $\al_1^{\mu_1}w_n\in W$. Since
$|k|\ge n+1\ge 2$, choose $0\ne \al_1\in k$. By assumption
$\al_1^{\mu_1}w_n \in W$. Since $W$ is a vector space, this
implies that $w_n\in W$ as claimed.

\medskip
Assume now that the assertion is true for $m-1$.
Write~\eqref{multi.v.d.m.1} as follows. Let $r=\mu_1$, then
\begin{equation}
\sum_{\mu_1+ \ldots +\mu_m =n}\al_1^{\mu_1}\cdots
\al_m^{\mu_m}w_{\mu_1,\ldots,\mu_m} =\sum_{r=0}^n\al_1^r\cdot \bar
w_r \in W,
\end{equation}
 where
\begin{equation}
\bar w_r =\sum_{\mu_2 + \ldots+\mu_m = n-r }\al_2^{\mu_2}\cdots
\al_m^{\mu_m}w_{r,\mu_2,\ldots,\mu_m}.
\end{equation}
By Lemma~\ref{v.d.m.2} all $\bar
w_r\in W$. Therefore by the induction assumption on $m-1$, all\\
$w_{r,\mu_1,\ldots,\mu_m}\in W$. Thus the assertion is true for
$m$, and the inductive step is complete.
\end{proof}

{\bf The proof of Proposition~\ref{v.d.m.1}} now follows from
Lemma~\ref{invar3} and from Proposition~\ref{main.vandermonde.1}.

\begin{cor}\label{2.18}
Let $A$ be a $k$-algebra with $|k|\ge n+1$ and let
$a_1,\ldots,a_m\in A$. Let $W \subseteq A$ be a subspace of $A$.
If
\begin{equation}\label{fold.1}
(\al_1 a_1+ \ldots + \al_m a_m)^n \in W
\end{equation}
for all $\al_1,\ldots,\al_m \in k$ then $P_n(a_1,\ldots,a_m)
\subseteq W$.
\end{cor}
\begin{proof}
Let
$W_x=W(x_1,\ldots,x_m)=\spn_k\{(\al_1x_1+ \ldots+\al_mx_m)^n\mid
\al_1,\ldots,\al_m\in k \}$, and similarly
$W_a=W(a_1,\ldots,a_m)=\spn_k\{(\al_1a_1+ \ldots+\al_ma_m)^n\mid
\al_1,\ldots,\al_m\in k \}$, so by assumption $W_a\subseteq W$.
Since $|k|\ge n+1$, by Proposition~\ref{v.d.m.1}
$W_x=P_n(x_1,\ldots,x_m) $. Substituting $x_r$ by $a_r$ implies
that $W_a=P_n(a_1,\ldots,a_m) $, which completes the proof.
\end{proof}

\begin{prop}\label{Sij-nil}
Let $A$ be an algebra over an infinite field $k$, and let
$a_1,\ldots,a_m \in A$. Then the subspace
$U=ka_1+ \ldots+ka_m\subseteq A$ is nil of bounded index $\le n$,
if and only if $P_n(a_1,\ldots,a_m)=\{0\}$.
\end{prop}

\begin{proof}
By definition $U=ka_1+ \ldots+ka_m$ is nil of bounded index $\le n$
if and only if for all $\al_1,\ldots,\al_m \in k$,
$(\al_1a_1+ \ldots+\al_m a_m)^n=0$. Thus the ``if" part of the
proposition follows from Lemma~\ref{invar3},
and the reverse assertion is a direct consequence of
Corollary~\ref{2.18}
with $W=\{0\}$.
\end{proof}

\section{Algebraicity of bounded degree}
Recall the space $P_n(x_1,\ldots,x_m)$ of Definition~\ref{DefSn}.
\begin{defn}
We denote \[P_{\le
r}(x_1,\ldots,x_m)=\sum_{n=1}^rP_n(x_1,\ldots,x_m)=\bigoplus
_{n=1}^rP_n(x_1,\ldots,x_m),\] and
\[P_{\ge
r}(x_1,\ldots,x_m)=\sum_{n=r}^{\infty}P_n(x_1,\ldots,x_m)=\bigoplus
_{n=r}^{\infty}P_n(x_1,\ldots,x_m).\]

Similarly, given a $k$-algebra $A$ and $a_1,\ldots,a_m\in A$, let\\
$P_{\le r}(a_1,\ldots,a_m)=\sum_{n=1}^rP_n(a_1,\ldots,a_m)$
and
$P_{\ge r}(a_1,\ldots,a_m)=\sum_{n=r}^{\infty}P_n(a_1,\ldots,a_m)$.
\end{defn}

\begin{lem}\label{Sij-alg2}
Let $D\ge d$. Let $A$ be a $k$-algebra with $|k|\ge D+1$ and let
$a_1,\ldots,a_m\in A$. Assume that for all $\al_1,\ldots,\al_m\in
k$, $\al_1 a_1+ \ldots +\al_ma_m$ is algebraic of degree at most
$d$.  Then $P_D(a_1,\ldots,a_m) \subseteq P_{\leq
d-1}(a_1,\ldots,a_m)$.

\end{lem}
\begin{proof}  Denote $P_j=P_j(a_1,\ldots,a_m)$ for any $j$, and denote $P_{\ge j}$ similarly. Let $u=\al_1a_1+ \ldots+\al_ma_m$. By
algebraicity $u^d$ is a linear combination of lower powers of $u$,
therefore $u^d\in P_{\le d-1}$. By Corollary~\ref{2.18} with
$W=P_{\le d-1}$,
it follows that $P_d\subseteq P_{\le d-1}$.

\medskip
If $D+1\ge d+2$, we can continue: By the same argument
$P_{d+1}\subseteq P_{\le d}$, while  the previous step
implies that $P_{\le d}\subseteq P_{\le d-1}$, hence
$P_{d+1}\subseteq P_{\le d-1}$. Continuing this way we
finally get that $P_D \subseteq P_{\leq
d-1}$.
\end{proof}

\begin{cor}\label{2.21}
Let $k$ be infinite, let $A$ be a $k$ algebra and let
$a_1,\ldots,a_m\in A$. If each element of $ka_1+ \ldots+ka_m$ is
algebraic of degree at most $d$ then\\ $P_{\ge 0}(a_1,\ldots,a_m)
\subseteq P_{\le d-1}(a_1,\ldots,a_m)$ and in particular,
$\dim_kP_{\ge 0}(a_1,\ldots,a_m) <\infty$.
\end{cor}

{ \begin{defn} Denote \[M_{d,m}={d+m-1\choose m}.\]
\end{defn}
 We note that
\begin{equation}\label{Sn-dim}
\dim_k P_n(x_1,\ldots,x_m) =
{{m+n-1}\choose {m-1}}
\end{equation}
 and
\begin{equation}\label{S-le-dim}
\dim_k P_{\leq n}(x_1,\ldots,x_m) = \sum^n_{j=0}\dim_k P_j =
\sum^n_{j=0}{j+m-1 \choose m-1}= {n+m\choose m}.
\end{equation}

Thus, given an algebra $A$ and elements $a_1,\ldots,a_m\in A$, we
have
\begin{equation}\label{S-le-dim2}
\dim_k P_{\leq d-1}(a_1,\ldots,a_m)\le M_{d,m}.
\end{equation}

\begin{lem}\label{Sij-alg4}
If $P_{\le M_{d,m}}(a_1,\ldots,a_m) \subseteq P_{\le
d-1}(a_1,\ldots,a_m)$ for some $d\ge1$ then the subspace $k
a_1+ \ldots+k a_m$ is algebraic of degree at most $M_{d,m}$.
\begin{proof}
Let $a=\al_1 a_1+ \ldots+\al_m a_m \in k a_1+ \ldots+k a_m$. Then by
assumption, for any $0\le n \le M_{d,m}$ we have (using
Lemma~\ref{invar3} with $x_i$ replaced by $a_i$)

\begin{equation}
a^n \in P_n(a_1,\ldots,a_m)\subseteq P_{\le M_{d,m}}(a_1,\ldots,a_m) \subseteq P_{\le d-1}(a_1,\ldots,a_m).
\end{equation}
Thus
\begin{equation}
\dim_k P_k \{a^n \mid 0\le n \le M_{d,m}\} \le \dim_k P_{\le
d-1}(a_1,\ldots,a_m) \le M_{d,m}.
\end{equation}
Therefore the set $\{a^n \mid  0 \le n \le M_{d,m}\}$ is linearly
dependent over $k$, and hence $a$ is algebraic of degree at most
$M_{d,m}$.
\end{proof}
\end{lem}
The following corollaries are immediate.
\begin{cor}\label{Sij-alg5}
If $P_{\ge 0}(a_1,\ldots,a_m) \subseteq P_{\le
d-1}(a_1,\ldots,a_m)$ for some $d\ge1$ then $k a_1+ \ldots+k a_m$
is algebraic of degree at most $M_{d,m}$.
\end{cor}

\begin{cor}\label{Sij-alg6}
If $\dim_k P_{\ge 0}(a_1,\ldots,a_m) < \infty$ then $k
a_1+ \ldots+k a_m$ is algebraic of bounded degree.
\end{cor}


We can now prove
\begin{prop}\label{2.9}
{\cite[Lemma 16(i)]{RTA}} Let $A$ be an algebra over an infinite
field $k$, and let $a_1,\ldots,a_m \in A$. Then the subspace
$ka_1+ \ldots+ka_m$ is algebraic of bounded degree if and only if
$\dim_k (P_{\ge 0}(a_1,\ldots,a_m))<\infty$.
\end{prop}

\begin{proof} One direction of the proof is given by
Corollary~\ref{2.21}, and the opposite direction -- by
Corollary~\ref{Sij-alg6}.
\end{proof}

\section{Filtered Algebraic Algebras}

Let $A$ be a filtered $k$-algebra with the filtration $\{F_n\}_{n\ge 0}$. We aim to prove that if $k$ is uncountable and $A$ is algebraic then every element of $(\gr A)_{\ge 1}$ is nilpotent. In particular, being graded-nil does imply being nil, for
 associated graded algebras of algebraic algebras over uncountable fields.

We first note that from (\ref{assoc-mult}) it follows (by induction) that if
$a_{p_1},\ldots,a_{p_n}$ satisfy\\ $a_{p_i}\in F_{p_i}$ for all
$1\le i \le n$ then
\begin{equation}\label{assoc-mult-more}
(a_{p_1}+F_{p_1-1})(a_{p_2}+F_{p_2-1})\cdots (a_{p_n}+F_{p_n-1})=
a_{p_1}a_{p_2}\cdots a_{p_n} + F_{p_1+\ldots + p_n - 1}.
\end{equation}

We note the following properties of the order-symmetric
polynomials in graded and filtered algebras. If the $k$-algebra $A$ is filtered by $\{F_n\}_{n\ge0}$
and $a_p,\ldots,a_q\in A$ satisfy $a_i\in F_i$ for all $p\le i\le q$, then
\begin{equation}\label{os-filt}
\s_{i_1,\ldots,i_m}(a_p,\ldots,a_q) \in F_{i_1 p+ \ldots+i_m q},
\end{equation}

where $m=q-p+1$. Similarly, suppose $B=\bigoplus_{i\ge 0}B_i$ is a graded $k$-algebra. If $b_p,\ldots,b_q$
are $m$ homogeneous elements with $b_i \in B_i$ then
\begin{equation}
\s_{i_1,\ldots,i_m}(b_p,\ldots,b_q) \in B_{i_1 p +\ldots+i_m q}.
\end{equation}

Now let $A$ be a $k$-algebra filtered by $\{F_n\}_{n\ge 0}$ and let
$B=\gr(A)$ be its associated graded algebra. Let $a_p,\ldots,a_q\in
A$ satisfy $a_i\in F_i$ for each $p\le i \le q$, and let
$b_i=a_i+F_{i-1}\in B$. We first note that by
(\ref{assoc-mult-more}), if $p\le p_j \le q$ for all $1\le j \le n$
then
\begin{equation}
b_{p_1}\cdots b_{p_n} = a_{p_1}\cdots a_{p_n} + F_{p_1+\ldots
+p_n-1}.
\end{equation}
Therefore if this product contains exactly $i_j$ occurrences of
$b_j$ for each\\ $p\le j \le q$, then
\begin{equation}
b_{p_1}\cdots b_{p_n} = a_{p_1}\cdots a_{p_n} + F_{i_p p+\ldots +i_q
q-1}.
\end{equation}
Since each $\s_{i_p,\ldots,i_q}(b_p,\ldots,b_q)$ is a sum of such
products, we have

\begin{equation}\label{multprob}
\s_{i_p,\ldots,i_q}(b_p,\ldots,b_q)=\s_{i_p,\ldots,i_q}(a_p,\ldots,a_q)+F_{pi_p+\ldots+qi_q-1}.
\end{equation}

\begin{thm}\label{my1}
Let $A$ be a filtered $k$-algebra, with a filtration
$\{F_n\}_{n\ge0}$. Let\\ $B=\gr(A)$ be its associated graded algebra.
If $A$ is LBD, then $B_{\ge1}$ is LBI. In particular, if $k$ is
uncountable and $A$ is a filtered algebraic $k$-algebra then
$\gr(A)_{\ge 1}$ is nil (and LBI).
\begin{proof}
Let $S$ be a finite-dimensional subspace of $B$. We may assume\\ $S=B_p \oplus
\ldots \oplus B_q$, with $1\le p \le q$, since every finite-dimensional subspace
of $B_{\ge 1}$ is contained in a subspace of this form. Let
$b=b_p+\ldots+b_q$ be any element of $S$, where $b_i \in
\dspfrac{F_i}{F_{i-1}}$ (we take $F_{-1}=\{0\}$). That is, for each
$i \in \{ p, \dots ,q \}$, we have $b_i=a_i+F_{i-1}$, with $a_i \in
F_i$. By assumption, the subspace $F_p+ \ldots +F_q$ is algebraic of
bounded degree, say $d$. Therefore by Corollary \ref{2.21}, we
have
\begin{equation}
\sps_{\ge0}(a_p,\ldots,a_q) \subseteq \sps_{\le
d-1}(a_p,\ldots,a_q).
\end{equation}
Now $\sps_{\le
d-1}(a_p,\ldots,a_q)$ is spanned by elements of the form $p_{i_1,\ldots,i_m}(a_p,\ldots,a_q)$, with $i_1+\ldots+i_m=d-1$.  Such $i_1,\ldots, i_m$ satisfy $i_1 p + \ldots + i_m q \le q(d-1).$
Therefore by (\ref{os-filt}) this implies
\begin{equation}
\sps_{\ge0}(a_p,\ldots,a_q) \subseteq F_{(d-1)q}.
\end{equation}

 Now let
\begin{equation}\label{N}
N= \Big\lceil{\frac{(d-1) q}{p}+1} \Big \rceil
\end{equation}
and suppose $i_p+\ldots+i_q\ge N$. Then
\begin{equation}
pi_p+\ldots+qi_q -1 \ge pN-1 \ge (d-1)q.
\end{equation}
Therefore
\begin{equation}
\s_{i_p,\ldots,i_q}(a_p,\ldots,a_q) \in F_{(d-1)q} \subseteq F_{p i_p + \ldots + q i_q-1},
\end{equation}
and so
\begin{equation}\label{eq0}
\s_{i_p,\ldots,i_q}(b_p,\ldots,b_q) =\s_{i_p,\ldots,i_q}(a_p,\ldots,a_q) + F_{p i_p + \ldots +q
i_q-1} =0_B.
\end{equation}

Thus $\sps_N(b_1,\ldots,b_n)=\{0_B\}$, and it follows from Proposition \ref{Sij-nil} that\\
$(b_p+\ldots+b_q)^N=0$. This shows that $S$ is nil of bounded index $N$. Thus $B$ is LBI.
The final statement follows from Theorem \ref{LBI-LBD.intro}.
\end{proof}
\end{thm}

\begin{remark}\label{rmk1}
We can use the bound $N$ of (\ref{N}) to derive a relation between the degrees of algebraicity in $A$ and the indices of nilpotence in $\gr(A)$. Namely,
if each $F_i$ is algebraic of bounded degree at most $d_i$, then each
$\dspfrac{F_1}{F_0} \oplus \dots \oplus \dspfrac{F_r}{F_{r-1}}$ is
nil of bounded index at most $(d_r-1) r + 1 \leq r  d_r$.
\end{remark}

\section{Integrality}
\begin{defn} For a $k$-algebra $A$, let $A[x]$ be the algebra of polynomials over $x$ with coefficients in $A$.
We say that an element $a(x)\in A[x]$ is \emph{integral} of degree
$n$ over $k[x]$ if $a(x)^n + p_{n-1}(x)a(x)^{n-1} + \ldots + p_0(x)
= 0$ for some $n\ge 1$ and polynomials $p_i(x)\in k[x]$.
\end{defn}
Note that if $a(x)$ is integral over $k[x]$ then in particular it is algebraic over $k(x)$, the field of rational functions over $k$.

\begin{lem}\label{thelem}
Let $A$ be a filtered $k$-algebra, with filtration $\{A_n\}_{n\ge0}$. Let \\$R=\bigoplus_{n\ge 0} x^n A_n  \subset A[x]$ be the Rees algebra.
Let $a(x)=a_1 x + ... + a_m x^m \in R$, where $a_i\in A_i$, and suppose that  $a(x)$ is integral over $k[x]$. Then $a(x)^N \in xR$ for some $N\ge 0$.
\begin{proof}
Suppose $a(x)$ is integral over $k[x]$ of degree $n$. Then for any $N\ge 0$
we have
\begin{equation}\label{inteq}
a(x)^N =  q_{n-1}(x)a(x)^{n-1} + \ldots + q_1(x)a(x) +q_0(x)
\end{equation}
for some polynomials $q_i(x)\in k[x]$. Choose $N\ge m(n-1)+1$.
The left-hand side of \eqref{inteq} is in $x^N A[x]$. The right-hand side is in $A_{m(n-1)}[x]$ (i.e., it is a polynomial in $x$ with coefficients in the subspace $A_{m(n-1)}$). Thus $a(x)^N \in x^N A[x] \cap A_{m(n-1)}[x] \\ \subseteq x^N A_{m(n-1)}[x] \subseteq x^N A_{N-1}[x] \subseteq \bigoplus_{n\ge 1} x^nA_{n-1} =xR$.
\end{proof}
\end{lem}
\begin{cor}\label{thecor}
Let $A$ be a filtered algebra and suppose that every element of the Rees algebra $R=\bigoplus_{n\ge 0} x^n A_n$ is integral over $k[x]$. Then $\gr(A)_{\ge1}$ is nil.
\begin{proof}
The map $\phi:\gr(A) \rightarrow  \frac{R}{xR} $ defined by\\
$\phi(a_0 + (a_1+A_0) + \ldots + (a_m+A_{m-1})) = a_0 + a_1 x + \ldots + a_m x^m + xR$\\
is a graded-algebra isomorphism. Under this map, an element\\ $(a_1+A_0)+\ldots+(a_m+A_{m-1})\in \gr(A)_{\ge 1}$ ($a_i\in A_i$) corresponds to\\ $a_1 x + \ldots + a_m x^m + xR \in \frac{R}{xR}$. By Lemma \ref{thelem} all such elements are nilpotent.
\end{proof}
\end{cor}
\begin{remark}
If $A$ is a filtered $k$-algebra which satisfies the LBD property (for example, if $k$ is uncountable and $A$ is algebraic over $k$), then by
\cite[Lemma 6]{AIF}, the extension algebra $A \otimes_k k(x)$ is LBD over $k(x)$.
\end{remark}

This raises the following question.
\begin{question}
If $A \otimes_k k(x)$ is LBD over $k(x)$, is every element of $A[x]$ integral over $k[x]$?
\end{question}
In view of Corollary \ref{thecor} (and since $R\subset A[x]$), a positive answer would give another proof of Theorem \ref{my1}.

\end{document}